\documentclass[runningheads]{llncs}

\usepackage[T1]{fontenc}
\usepackage{newtxtext}
\usepackage[varvw]{newtxmath}
\usepackage{amsmath}

\usepackage{enumitem}
\usepackage{bbding}
\usepackage{graphicx}
\usepackage[colorlinks,linkcolor=blue,citecolor=blue,urlcolor=blue,hypertexnames=true]{hyperref}
\usepackage{color}

\usepackage[margin=1.5in]{geometry}

\newcommand{\f}{\mathcal{F}}

\newcommand{\lang}{\mathcal{L}}
\newcommand{\logl}{\mathscr{L}}
\newcommand{\bvdash}{\mathrel{\scalebox{1.2}{$\vdash$}}}
\newcommand{\nb}{\neg_B}
\newcommand{\nd}{\neg_D}
\newcommand{\sr}{\mathscr{S}}
\newcommand{\mm}{\mathcal{M}}
\newcommand{\e}{\mathbf{E}}
\newcommand{\cb}{\circ_B}
\newcommand{\hstar}{\mathord{\raisebox{0.08ex}{\scalebox{0.65}{\textup{\FiveStarOpen}}}}}
\newcommand{\hstarB}{\mathord{\raisebox{0.08ex}{\scalebox{0.65}{\textup{\FiveStarOpen}}}}_{\!B}}
\newcommand{\db}{\hstarB}
\newcommand{\hstarD}{\mathord{\raisebox{0.08ex}{\scalebox{0.65}{\textup{\FiveStarOpen}}}}_{\!D}}
\newcommand{\cd}{\circ_D}
\newcommand{\dD}{\hstarD}

\begin{document}
\title{From Non-normality to Paranormality: Modal Negations in the Logic \texorpdfstring{$\e$}{E}}
\author{Esha Jain\orcidID{0009-0009-5648-030X} }
\authorrunning{E. Jain}
%
\institute{Department of Mathematics,
  Indraprastha Institute of Information Technology-Delhi\\
  New Delhi, India\\
\email{eshaj@iiitd.ac.in}
}
\maketitle
\begin{abstract}
In this article, we study two negations induced by modal operators in $\e$, the smallest classical modal logic, namely $\neg \Box$ and $\neg \Diamond$. Using neighbourhood semantics, we prove that these negations satisfy the paraconsistency and paracompleteness criteria proposed by both Marcos and Omori. In fact, with respect to these negations, $\e$ belongs to a special class of paraconsistent logics called logics of formal inconsistency and its dual class, logics of formal undeterminedness (a class of paracomplete logics). We explore several other properties of these modal negations. We examine how the behaviour of these modal negations changes when $\e$ is strengthened and compare $\e$ with a few familiar non-classical logics.

\end{abstract}
\section{Introduction}

Modal logics are primarily used to qualify the truth of a proposition. This is done through notions such as necessity and possibility, knowledge and belief, permission and obligation, etc. Conventionally, necessity and possibility are denoted by $\Box$ and $\Diamond$, respectively. In the present scenario, modal logic has applications in various fields such as computer science, artificial intelligence, multi-agent systems, and formal philosophy (see, for example, \cite{Pnueli1977,RaoGeorgeff1995,vonWright1951}).
While applying modal logics in practical contexts, logicians aim to remain as rational and realistic as possible. Despite having a strong mathematical structure, classical propositional logic ($\mathbf{CPL}$) fails to capture many real-world aspects, such as handling incomplete or contradictory information. This limitation gives rise to the need to investigate paraconsistent and paracomplete logics through modal frameworks. 

A logic $\mathscr{L}=\langle\lang,\bvdash_{\mathscr{L}}\rangle$, where $\lang$ is the set of formulas and $\bvdash_\logl$ is the consequence relation of the logic, is said to be \emph{paraconsistent} if it invalidates a \emph{principle of logical explosion}. The most common principle of logical explosion is expressed as follows: for all $\alpha,\beta\in\lang$, $\{\alpha,\neg\alpha\}\bvdash_{\mathscr{L}}\beta$. This is called \emph{ECQ (ex contradictione sequitor quodlibet)}.  There are other principles of explosion (see \cite{MR4939722}) that produce other kinds of paraconsistency. The dual of a paraconsistent logic is called a paracomplete logic.  A logic $\mathscr{L}=\langle\lang,\bvdash_\logl\rangle$ is said to be \emph{paracomplete} if the \emph{law of excluded middle (LEM)} fails in it, i.e., there exists $\alpha\in\lang$ such that $\not\bvdash_\mathscr{L}\alpha\lor\neg\alpha$. In a paracomplete logic, two propositions $\alpha$ and $\neg\alpha$ can both be false.

To the best of our knowledge, Ja\'{s}kowski (see \cite{jaskowski1948rachunek,jaskowski1969propositional}) first introduced a formal paraconsistent logic, known as Discussive (or Discursive) logic $\mathbf{D_2}$, using the normal modal logic $\mathbf{S5}$. He replaced material implication and classical conjunction with two connectives based on the modal operator of possibility $\Diamond$.
In subsequent foundational work \cite{beziau2002s5}, B\'{e}ziau claimed that $\mathbf{S5}$ is paraconsistent while considering the modal negation possibly not ($\Diamond\neg$). In \cite{dosen1984negative}, Do\v{s}en studied negations induced by modal operators corresponding to impossibility and non-necessity over a fragment of intuitionistic logic.  Marcos  \cite{MR2193566} expanded B\'{e}ziau's framework by defining a criterion in semantic terms to determine paraconsistent and paracomplete behaviour of general operators or connectives in the language. Recently, Omori \cite{homori2025s5} reconsidered the question of when a modal logic should be regarded as paraconsistent or paracomplete, emphasizing that paraconsistency and paracompleteness should be examined with respect to both possibly not ($\Diamond\neg$) and necessarily not ($\Box\neg$) forms of negation. However, a common feature in most of these works is that they are carried out in the setting of normal modal logics (i.e. systems extending $\mathbf{CPL}$ with the (K) axiom ($\Box(p\rightarrow q)\rightarrow (\Box p \rightarrow \Box q)$) and the rule of necessitation \ref{nec rule}). Little attention has been paid to non-normal logics. The need for non-normal modal logic arises when the closure principles of normal modal logics are too strong for the intended interpretation of the modality. 

While paraconsistency and paracompleteness are commonly characterized by the failure of ECQ and LEM, respectively, these notions are inadequate for determining the behaviour of operators induced by modalities. These definitions do not establish whether a defined unary operator acts as a negation. To overcome this, in the present work, we utilize the conceptual frameworks introduced by Marcos and Omori in \cite{MR2193566,homori2025s5}  to investigate the behaviour of negations induced by modal operators within $\e$, a non-normal modal logic,  also called the smallest classical modal logic (Definition \ref{def:normal  and minimallogic}). We establish that these induced negations exhibit paraconsistent and paracomplete behaviour within $\e$. We then show that $\e$ belongs to a special class of paraconsistent logics called logics of formal inconsistency (LFIs) (Definition \ref{def:lfi}) introduced in \cite{MR2008227}. In LFIs,  the notion of consistency is formalised in the object language. The dual of LFI is called LFU, logics of formal undeterminedness (Definition \ref{def:lfu}). We establish that, despite weak structure,  $\e$ is both an LFI and an LFU with respect to modal negations.  Finally, we compare the non-classical behaviour of $\e$ with a few well-known logics. 
 
The article is structured as follows. We describe the logic $\e$ along with the required definitions and results in Section \ref{sec:system e}. Then, we present the notions introduced by Marcos and Omori in Section \ref{sec:paraconsis and pcompleteness} along with the results on the operators induced by the modality in $\e$. In Section \ref{sec:lfi and lfu}, we prove that $\e$ gives rise to two LFIs and two LFUs. Finally, we compare the paraconsistent and paracomplete behaviour of $\e$ with other well-known logics in Section \ref{sec:comparison}.

\section{The logic $\e$}\label{sec:system e}
\textbf{Language.}
\label{language}Let  $\mathrm{Var}=\{p,q, p_1,\dots\}$ be a denumerable set of propositional variables and $\{\land,\neg, \Box\}$ be a set of connectives. Let $\mathcal{L'}$ be the language consisting of all the formulas generated over $\mathrm{Var}$ using $\{\land,\neg, \Box\}$. 
Abbreviations for the connectives $\lor$ (disjunction), $\rightarrow$ (implication) and $\leftrightarrow$ (equivalence) are standard. 
 Since these operators can be defined in terms of $\land$ and $\neg$, the formal definition of the language uses only these three connectives. For any formula $\alpha$, the dual of an operator $\Box$ is defined by $\Diamond \alpha:=\neg \Box \neg \alpha$.

 To classify systems of modal logics, we first present a few formulas, rules, and definitions. 

 \begin{enumerate}[label=(\roman*)]
     \item \label{k axiom} $\Box(\alpha\rightarrow\beta)\rightarrow (\Box \alpha\rightarrow \Box \beta)$\hfill  (K)
     \item \label{M axiom}$\Box(\alpha\land\beta)\rightarrow (\Box\alpha\land \Box \beta)$\hfill (M)
     \item \label{C axiom}$(\Box\alpha\land \Box \beta)\rightarrow\Box(\alpha\land\beta) $\hfill(C)
     \item \label{nec rule}$\dfrac{\alpha}{\Box\alpha}\;$ if $\alpha$ is a theorem\hfill (Nec)
     \item \label{RE rule}$\dfrac{\alpha\leftrightarrow\beta}{\Box \alpha\leftrightarrow\Box \beta}\; $ if $\alpha\leftrightarrow\beta$ is a theorem; \hfill (RE)
    \item \label{MP rule}$\dfrac{\alpha\;\;\alpha\rightarrow\beta}{\beta}$\hfill (MP)
     \end{enumerate}
 The following definition is adapted from \cite{pacuit2017neighborhood}.
 \begin{definition}\label{def:normal  and minimallogic}
    A logic $\mathscr{L}=\langle\mathcal{L'},\bvdash_\logl\rangle$ is called a \emph{normal modal logic} if all instances of classical propositional tautologies and all instances of (K) are theorems of $\sr$ and it is closed under the rules (Nec) and (MP).

      A logic $\mathscr{L}=\langle\mathcal{L}',\bvdash_\logl\rangle$  is called a \emph{minimal modal logic} (or classical modal logic in \cite{chellas1980modal}) if all instances of classical propositional tautologies are theorems of $\sr$ and it is closed under the rules (RE)  and (MP).
 \end{definition}

  The well-known logic $\mathbf{K}$ is obtained by adding the axiom scheme (K) and the rule (Nec) to $\mathbf{CPL}$. It is the smallest normal modal logic (in terms of the number of theorems). 
 
 Non-normal modal logics have a variety of applications in epistemic, deontic and evidential reasoning, e.g., in epistemic logic, if we consider the smallest normal modal logic $\mathbf{K}$ and interpret $\Box$ as `the agent knows that', then formulas like $\Box p\land \Box q\rightarrow \Box(p \land q)$ become valid. This shows that agents are perfect reasoners in $\mathbf{K}$.  Normal modal logics suffer from the problem of logical omniscience, which makes it difficult to reason about agents that are not perfect reasoners. To study this in detail, one can see \cite{pacuit2017neighborhood,holliday2015epistemic}. 

In this article, we focus on the non-normal\footnote{Observation 2.32 in Chapter 2 of \cite{pacuit2017neighborhood} shows that (K) is not valid in $\e$.} modal system $\e$.

 The logic $\e=\langle\mathcal{L'}, \bvdash_\e\rangle$ is defined by the following Hilbert calculus.

 \noindent\textbf{Axiom Schema.}
  All instances of tautologies of classical propositional logic ($\mathbf{CPL}$).
    
    \noindent
    \textbf{Inference Rules.} (RE) and (MP).
    
\begin{definition}\label{def:syntactic consequence e}
     Let $\Gamma\cup\{\alpha\}\subseteq\lang'$. A \emph{derivation} of $\alpha$ from $\Gamma$ in $\e$ is a finite sequence $(\alpha_{1},\ldots,\alpha_{n})$ of elements in $\lang'$, where $\alpha_{n}=\alpha$ and for each $1\le i\le n$, $\alpha_{i}$ is an instance of an axiom of $\e$, or a member of 
     $ \Gamma$, or obtained  by an application of (MP) on $\alpha_{j},\alpha_{k}$ 
     where $1\le j,k<i$, or obtained by an application of the rule (RE) on $\alpha_j$ where $j<i$, and $\alpha_j$ is a theorem.

We say that $\alpha$ is \emph{syntactically derivable} or \emph{syntactically entailed} from $\Gamma$, and write $\Gamma\bvdash_\e\alpha$, if there is a derivation of $\alpha$ from $\Gamma$. If $\Gamma=\emptyset$, then $\alpha$ is called a theorem.
     \end{definition}
    
\begin{definition}\label{def:structural logic}
 Let $\logl=\langle\lang,\bvdash_\logl\rangle$ be a logic whose language $\lang$ is generated over $\mathrm{Var}$. A \emph{substitution} is defined as any function $\sigma:\mathrm{Var}\to\lang$ that extends to a unique endomorphism (also denoted by $\sigma$) from $\lang$ to itself via the universal mapping property. 
    The logic $\logl$ is called \emph{structural} if for every $\Gamma\cup\{\alpha\}\subseteq\lang$, $\Gamma\bvdash_\logl\alpha$ implies $\sigma[\Gamma]\bvdash_\logl\sigma(\alpha)$, for every substitution $\sigma$ of formulas for variables.
\end{definition}
     \begin{remark}\label{rem:e structural} The logic $\e$ is structural. This follows by induction on the length of a derivation.
     According to Definition \ref{def:normal  and minimallogic}, $\e$ is a minimal modal logic. In fact, $\e$ is the smallest classical modal logic.  $\e$ is sound and complete with respect to the class of all neighbourhood frames \cite[Theorem 2.61] {pacuit2017neighborhood}. In this article, we 
     prove most of the results using neighbourhood  semantics.  So, we present a few notions of neighbourhood  semantics. 
     \end{remark}

\noindent
\textbf{Neighbourhood semantics} were developed independently by Scott and Montague \cite{scott1970advice,montague1970universal}. Since then,
 these semantics have been extensively studied (see, for example,  \cite{segerberg1971essay,chellas1980modal,hansen2009neighbourhood,pacuit2017neighborhood}). The definitions provided in this section are adapted from \cite{pacuit2017neighborhood}.  

\begin{definition}
    Let $W$ be a non-empty set. A function $N: W\rightarrow\mathcal{P}(\mathcal{P}(W))$ is called a \emph{neighbourhood function}. A pair $\langle W, N\rangle$ is called a \emph{neighbourhood frame} if $W$ is a non-empty set and $N$ is a neighbourhood function. In modal logic, elements of $W$ are called worlds or states.
\end{definition}

 \begin{definition}
     Suppose that $\mathcal{F}=\langle W,N\rangle$ is a neighbourhood frame. A \emph{neighbourhood model} based on $\mathcal{F}$ is a tuple $\langle W,N, V\rangle$, where $V:\mathrm{Var}\rightarrow \mathcal{P}(W)$ is a valuation function (a function assigning a set of states to each propositional variable).
\end{definition}
\begin{definition}\label{def:truth }
     Suppose that $\mathcal{M}=\langle W,N,V\rangle $ is a neighbourhood model, and $w\in W$. Suppose $\alpha, \beta\in \mathcal{L'}$. Then, \emph{truth} of a formula in the model at world $w$ is defined as follows.
    \begin{enumerate}[label=(\roman*)]
        \item $\mathcal{M},w\models p$ iff $w\in V(p)$ where $p\in \mathrm{Var}$;
        \item $\mathcal{M}, w\models \neg \alpha$ iff $\mathcal{M},w\not\models \alpha $;
        \item $\mathcal{M}, w\models \alpha\land \beta$ iff $\mathcal{M},w\models \alpha $ and $\mathcal{M},w\models \beta $;
        \item $\mathcal{M},w\models \Box\alpha $ iff $\llbracket\alpha\rrbracket_\mathcal{M}\in N(w)$.

    \end{enumerate}  
Here, $\llbracket \alpha\rrbracket_\mathcal{M}$ denotes the truth set of $\alpha$ in model $\mm$, i.e., $\llbracket\alpha\rrbracket_\mathcal{M}=\{s\in W:\mathcal{M},s\models\alpha \}$. 
\end{definition}

For any $X\subseteq W$, let $X^C$ denote the complement of $X$, i.e., $X^C:=W\setminus X$.
Using the above definition of truth, we have \label{truth of diamond operator}  $\mathcal{M},w\models \Diamond\alpha $ iff $\llbracket\alpha\rrbracket^C_\mathcal{M}\not\in N(w)$.

\begin{definition}
    Let $\mathscr{N}$ be the class of all neighbourhood  models. Let $\Gamma\cup \{\alpha\}\subseteq \lang '$. Then, $\Gamma$ \emph{semantically entails} $\alpha$, denoted by $\Gamma\models_\e \alpha$, if for every model $\mm\in \mathscr{N}$ and every world $w$ of $\mm$, $$\mm,w\models \beta \hbox{ for every } \beta\in \Gamma \hbox{ implies } \mm, w\models \alpha.$$ 
    If $\Gamma=\emptyset$, then we write $\models\alpha$; equivalently, $\alpha$ is called valid in $\e$ if $\mm,w\models \alpha$  for every   $\mm\in \mathscr{N}$ and every world $w$ of $\mm$.
\end{definition}

 In this article, we use $\models$ for truth at a world in a model and $\models_\e$ for semantic consequence in $\e$. The subscript $\e$ is omitted from $\bvdash_\e$ and $\models_\e$ whenever the underlying logic is clear.
 
\noindent Now, in the smallest classical modal logic $\e$, we define two unary operators. For any formula $\alpha\in \mathcal{L'}$, 
$$\neg_B\alpha:=\neg\Box\alpha\hbox{ and }\nd\alpha:=\neg\Diamond\alpha.$$

The subscripts $B$ and $D$ serve as natural mnemonics for the modalities $\Box$ and $\Diamond$, respectively. $\nb$ represents classical negation applied to the modal operator $\Box$, and $\nd$ represents classical negation applied to the modal operator $\Diamond$. In
\cite{MR2193566}, $\nb$ and $\nd$ are expressed by $\smile$ and $\frown$ respectively.

The semantic conditions for $\nb$ and $\nd$ in $\e$ according to Definition \ref{def:truth } are as follows. 
\begin{enumerate}[label=(\roman*)]
\item $\mm, w\models\nb\alpha$ iff $\llbracket\alpha\rrbracket\not\in N(w)$

\item $\mm, w\models\nd\alpha$ iff $\llbracket\alpha\rrbracket^C\in N(w)$.
\end{enumerate}

\section{Paraconsistency and Paracompleteness}\label{sec:paraconsis and pcompleteness}

A brief idea about paraconsistency and paracompleteness is provided in the introduction. In this section, we present formal definitions of the notions given by Marcos and Omori in \cite{MR2193566} and \cite{homori2025s5}, respectively, and explore properties of $\e$.

Given two logics $\sr_1=\langle\lang_1, \bvdash_1\rangle$ and $\sr_2=\langle\lang_2, \bvdash_2\rangle$,  $\sr_2$ is an \textit{extension} of  $\sr_1$ if $\lang_1\subseteq \lang_2$ and for any $\Gamma\cup \{\alpha\}\subseteq \lang_1$, $\Gamma\bvdash_1\alpha$ implies $\Gamma\bvdash_2 \alpha$\footnote{$\sr_2$ is a \textit{conservative extension} of  $\sr_1$ if $\lang_1\subseteq \lang_2$ and for any $\Gamma\cup \{\alpha\}\subseteq \lang_1$, $\Gamma\bvdash_1\alpha$ iff $\Gamma\bvdash_2 \alpha$.}.

 \textbf{Marcos's Perspective.} Marcos provided model-theoretic definitions for classifying a logic as paraconsistent, paracomplete or paranormal logic. The concept of decent negation explained below follows Marcos's multi-conclusion framework. The remainder of this section uses the standard single-conclusion consequence relation.  While Marcos formulated these definitions in the context of Kripke semantics, the definitions solely depend on a satisfaction relation. So, we adapt these definitions directly for neighbourhood  models of $\e$ without altering the model-theoretic definitions given in \cite{MR2193566}.

 Let $\sr=\langle \mathcal{L}, \models_{\mathscr{S}}\rangle$ be an arbitrary logic defined over the language $\mathcal{L}$, where $\models_{\mathscr{S}}$ is its semantic consequence relation.
 The multi-conclusion relation is defined as follows. For $\Gamma, \Delta\subseteq \lang$, 
 $\Gamma\models_\sr \Delta$ holds iff for every model $\mm$ and state $w$, $$\mm,w \models \gamma \hbox{ for every } \gamma\in \Gamma \hbox{ implies } \mm,w \models \delta \hbox{ for some } \delta\in \Delta.$$  
 
 Suppose $\sr$ is equipped with a unary operator $\sim$.  We define $\sim^0 \alpha:=\alpha$ and  $\sim^{n+1}\alpha:=\sim \sim^n \alpha$ for any formula $\alpha\in \lang$ and $n\in \mathbb{N}$. Then,  $\sim$ is called a \textit{decent negation}\label{dfn:decent negation} if it satisfies the following conditions.  
 
 For every $n\in \mathbb{N}$:
 \begin{itemize}[label=\textbullet]
     \item there exist $\Gamma,\Delta\subseteq \mathcal{L}$ and  $\alpha\in \mathcal{L}$ such that\\
$\Gamma,\sim^{n+1}\alpha\not\models_\sr\;\sim^n\alpha,\Delta$\hfill (\textit{n-verificatio})
     \item there exist $\Gamma,\Delta\subseteq \mathcal{L}$ and $\alpha\in \mathcal{L}$ such that\\  $\Gamma,\sim^{n}\alpha\not\models_\sr\;\sim^{n+1}\alpha,\Delta$\hfill (\textit{n-falsificatio})
 \end{itemize}

Based on this, Marcos classifies the given logic using the following properties. 
 
\begin{enumerate}[label=(\roman*)]
 \item  The logic $\sr$ is called \textit{overcomplete} if all of its models satisfy either of the following conditions. 
 \begin{enumerate}[label=(\alph*)]
     \item  For all $\alpha\in \lang$ and for all $w\in W$, $\mm,w\models\alpha$.
     \item  For all $\alpha\in \lang$ and for all $w\in W$, $\mm,w\not\models\alpha$.
 \end{enumerate}
 
\item The logic $\sr$ is said to be \textit{non-overcomplete} if there exist $\alpha$ and $\beta\in \lang$ such that $\alpha\not\models_\sr \beta$.
\item A model $\mm$ in  $\sr$ is said to be \textit{$\sim$-inconsistent} if there exists  $\alpha\in \lang$ and a $w\in W$ such that $\mm,w \models \alpha$ and $\mm,w \models \sim\alpha$.
\item A model $\mm$ in  $\sr$ is said to be $\sim$-\textit{incomplete} or $\sim$-\textit{undetermined} if there exists $\alpha\in \lang$ and a $w\in W$ such that $\mm,w \not\models \alpha$ and $\mm,w \not\models \sim\alpha$.

\end{enumerate}
  The logic $\sr$ is said to be \textit{negation-inconsistent} if some formula and its negation are both satisfied at the same world in some model of $\sr$. 
 
 The logic $\sr$ is said to be \textit{negation-incomplete} if some formula and its negation are both unsatisfied at the same world in some model of $\sr$.

\begin{definition}\label{def:marcos}
Suppose $\sr=\langle \mathcal{L}, \models_{\mathscr{S}}\rangle$ is a logic equipped with a decent negation. Then, 
\begin{enumerate}[label=(\roman*)]\item\label{def:pconsistent Marcos}
  $\sr$ is \emph{paraconsistent} if it is negation inconsistent yet non-overcomplete. 
  
 \item \label{def:pcomplete Marcos}  $\sr$ is \emph{paracomplete} if it is negation incomplete yet non-overcomplete.
  
\item \label{def:pnormal Marcos}  $\sr$ is \emph{paranormal} if it is both paraconsistent and paracomplete.  
\end{enumerate}
\end{definition}

\textbf{Omori's Perspective.} Omori proposed notions of  paraconsistency and paracompleteness for logics that have modalities and  are extensions of classical propositional logic (\textbf{CPL}), based on the failure of certain formulas involving modal negations. 

     \begin{definition}\label{def:omori} Let $\lang$ be a language generated over $\mathrm{Var}$ with a set of connectives $\{\neg,\land,\lor, \rightarrow, \Box,\Diamond\}$ (primitive or defined). Suppose $\sr=\langle \mathcal{L}, \bvdash_{\mathscr{S}}\rangle$ is a logic that is an extension of classical propositional logic $(\mathbf{CPL})$ and is closed under uniform substitution of arbitrary formulas for propositional variables. Then, \begin{enumerate} [label=(\roman*)]\item \label{def:pconsistent Omori}$\sr$ is said to be \emph{paraconsistent} iff   $\not\bvdash_\sr p\rightarrow (\Diamond\neg p
\rightarrow q) $ and  $\not\bvdash_\sr p\rightarrow (\Box\neg p
\rightarrow q) $ where $p,q\in \mathrm{Var}$.
\item \label{def:pcomplete Omori}
     $\sr$ is said to be \emph{paracomplete} iff  $\not\bvdash_\sr p\lor (\Diamond\neg p
) $ and $\not\bvdash_\sr p \lor \Box \neg p$ where $p\in \mathrm{Var}$.
\end{enumerate}
\end{definition}

\begin{remark}
    Two formulas $\alpha$, $\beta\in \lang'$ are equivalent in $\e$, denoted by $\alpha\equiv\beta$, iff $\alpha\models\beta$ and $\beta\models\alpha$. In $\e$, by the definition of $\Diamond$  in $\lang'$, available axioms and rules yield $\Diamond\neg\alpha\equiv \nb \alpha$. Similarly, $\Box\neg \alpha\equiv\nd\alpha$. 
\end{remark}

\textbf{Properties of modal negations in} $\e$. The following theorems present the behaviour of $\nb$ and $\nd$ in $\e$.

\begin{theorem}\label{thm:ecq fails} ECQ fails with respect to $\nb$ and $\nd$ in $\e$, i.e.,
there exist $\alpha$, $\beta\in\mathcal{L'}$, such that the following hold.

 \noindent
 (i) \label{ecq fails nb}$\alpha,\nb\alpha\not\models \beta$;\qquad
(ii)  \label{ecq fails nd}$\alpha,\nd\alpha\not\models \beta$.

\end{theorem}
\begin{proof}Let $W=\{w,s,t\}$ and $p$, $q$ be propositional variables, i.e., $p$, $q\in \mathrm{Var}$.  Suppose $\mm=\langle W,N,V \rangle$ where $N(w)=\{\{w,s\}, \{s,t\}\}$ and  
 $V$ is a valuation such that $V(p)=\{w\}$ and $V(q)=\{s\}$.  Then, we obtain  
    $\mm,w\models p$ , $\mm,w\models \nb p$ and $\mm,w\models \nd p$ but $\mm,w\not\models q$ since $w\not\in V(q)$. Therefore,  $\{p,\nb p\}\not\models q$ and  $\{p,\nd p\}\not\models q$.
\qed
\end{proof}

\begin{theorem}\label{thm:lem fails}
    Law of excluded middle (LEM) fails for both $\nb$ and $\nd$ in $\e$, i.e., there exists a formula $\alpha\in \mathcal{L'}$ such that 

    \noindent
    (i) $\not\models \alpha\lor \nb \alpha$;\qquad  (ii) $\not\models \alpha\lor \nd \alpha$.
    
\end{theorem}
\begin{proof}
Let $W=\{w,s,t\}$ and $p$ be a propositional variable, i.e., $p\in \mathrm{Var}$. Suppose $\mm=\langle W,N,V \rangle$ where $N(w)=\{\{w,s\}, \{s,t\}\}$ and  
 $V$ is a valuation such that $V(p)=\{s,t\}$.  Then, we obtain  
    $\mm,w\not\models p$ , $\mm,w\not\models \nb p$ and $\mm,w\not\models \nd p$ . Therefore,  $\not\models p \lor\nb p$ and  $\not\models p\lor \nd p$. 
    \qed
\end{proof}

\noindent
In the following lemma, we prove that $\nb$  and $\nd$ are decent negations in $\e$. 
\begin{lemma}\label{lem:decent neg}
$\nb$ and $\nd$ defined in $\e$ are decent negations.
\end{lemma}
\begin{proof}
    Suppose $W=\{w\}$ and $N(w)=\{\{w\}\}$. Let $\alpha\in\mathcal{L}'$ be any formula and $\f=\langle W, N\rangle$ be a frame.

    \noindent
   For any model $\mm$ based on $\f$,
$\llbracket\alpha\rrbracket_\mm=W$ if $\mm, w\models \alpha$ and $\llbracket\alpha\rrbracket_\mm=\emptyset$ if $\mm, w\not\models \alpha$.

We know that $\mm, w\models\Box \alpha $  iff $ \llbracket\alpha\rrbracket_\mm\in N(w)$. Then, it is easy to verify that $\llbracket\alpha\rrbracket_\mm\in N(w)$ iff $ \llbracket\alpha\rrbracket_\mm=W$. Also, we have $\llbracket\alpha\rrbracket_\mm=W$ iff $ \mm, w\models \alpha$.

Combining these, we obtain $\mm, w\models \Box \alpha $ iff $\mm,  w\models \alpha$.

\noindent
Since $\nb\alpha=\neg\Box\alpha$ and $\nd\alpha=\neg\Diamond\alpha$, using $\mm,w\models \Box \alpha$ iff $\mm, w\models \alpha$, we obtain $\mm, w\models \neg\Box \alpha$ iff $\mm, w\not\models \Box \alpha$ iff $\mm, w \not\models \alpha$ iff $\mm, w \models \neg \alpha$.

Similarly, we obtain $\mm, w\models \Box\neg \alpha$ iff $\mm, w\models \neg \alpha$ iff $\mm, w \not\models \alpha$.

\noindent
We note that both $\nb$ and $\nd$ behave like classical negation in the frame $\f=\langle W, N \rangle$.

\noindent
 We know that  $\neg^n\alpha\equiv\alpha$ if n is even and $\neg^n\alpha\equiv\neg \alpha$ if n is odd because $\neg$ is a classical negation.

\noindent
Let $p\in \mathrm{Var}$. Thus, for any model $\mm$ based on $\f$ and $\sim\in \{\nb,\nd\}$, we obtain
\begin{center}
$\mm, w\models \sim^n p  $ if and only if $\begin{cases}
   \mm, w\models p & n \hbox
    { is }\mathrm{even};\\
\mm, w\not\models p& n \hbox
    { is }\mathrm{odd}.    
\end{cases}$
\end{center}

\noindent
To show that $\nb$ and $\nd$ satisfy $n$-verificatio and $n$-falsificatio, we need to prove that $\sim^{n+1}\alpha\not\models \sim^n \alpha$,  and $\sim^{n}\alpha\not\models \sim^{n+1} \alpha$ for every $n$, where $\sim\in \{\nb,\nd\}$ (we take $\Gamma=\Delta=\emptyset$ in the definition of decent negation).

\noindent
Suppose  $\mm'=\langle W, N, V\rangle$ where $V(p)=\{w\}$. Then, $\mm',w\models p$ and hence, $\mm',w\not\models \neg p$. This implies $\mm', w\not\models \nb p$ and $\mm',w\not\models \nd p$.

Thus, if $n$ is odd, $\mm',w \models \sim^{n+1} p$ but $\mm',w \not\models \sim^n p$ for  $\sim\in \{\nb,\nd\}$. This proves $n$-verificatio for $\nb$ and $\nd$ if $n$ is odd.   If $n$ is even, we obtain $\mm',w \models \sim^n p$ but $\mm',w \not\models \sim^{n+1} p$ for $\sim\in \{\nb,\nd\}$.
This proves $n$-falsificatio for $\nb$ and $\nd$ if $n$ is even.

 \noindent
For the remaining cases, we take $\alpha$ as $\neg p$.  Since $\neg p$ is not true at $w$ in $\mm'$, the same argument establishes $n$-verificatio for even $n$ and  $n$-falsificatio for odd $n$ for both $\nb$ and $\nd$.   
  \qed
\end{proof}

 We have proved that $\nb$ and $\nd$ are decent negations. It is clear that $\e$ is an extension of $\mathbf{CPL}$ and is closed under uniform substitution, as it is structural (by Remark \ref{rem:e structural}). We are therefore in a position to determine the behaviour of $\e$ with respect to these negations according to the frameworks of Marcos and Omori. 
 
\begin{corollary}\begin{enumerate}[label=(\roman*)]
\item $\e$ is paraconsistent with respect to $\nb$ and $\nd$ in Marcos's view (Def. \ref{def:marcos}). 
\item $\e$  is paraconsistent in  Omori's view (Def. \ref{def:omori}).
\end{enumerate}
\end{corollary}
\begin{proof}Suppose $p,q\in \mathrm{Var}$. Suppose $\mm=\langle W,N,V \rangle$ where $W=\{w,s,t\}$, $N(w)=\{\{w,s\}, \{s,t\}\}$ and  
 $V$ is a valuation such that $V(p)=\{w\}$ and $V(q)=\{s\}$.  Then, we obtain  
    $\mm,w\models p$ , $\mm,w\models \nb p$ and $\mm,w\models \nd p$ but $\mm,w\not\models q$ since $w\not\in V(q)$. Therefore,  $\{p,\nb p\}\not\models q$ and  $\{p,\nd p\}\not\models q$.
 \begin{enumerate}[label=(\roman*)]\item  $\nb$ and $\nd$ qualify as  decent negations by Lemma \ref{lem:decent neg}. The model $\mm$ shows that $\e$ is $\nb$-inconsistent and $\nd$-inconsistent. Since $\mm,w\models p$, this implies $\mm,w \not\models \neg p$. Thus, $p \not\models \neg p$. Hence, $\e$ is non-overcomplete as well. Using the Def. \ref{def:marcos} \ref{def:pconsistent Marcos}, $\e$ is paraconsistent with respect to $\nb$ and $\nd$.  

\item Since $\mm, w\not\models p\rightarrow (\nb p\rightarrow q)$, we have  $\not\models p\rightarrow (\nb p\rightarrow q)$. Also, $\mm,w \not\models p\rightarrow (\nd p
\rightarrow q)$. Thus, we obtain $\not\models p\rightarrow (\nd p
\rightarrow q)$. By soundness of $\e$ with respect to neighbourhood semantics, $\e$ is paraconsistent according to Def. \ref{def:omori} \ref{def:pconsistent Omori}.  \qed 
\end{enumerate}

\end{proof}
\begin{corollary}
   \begin{enumerate}[label=(\roman*)]
\item $\e$ is paracomplete with respect to $\nb$ and $\nd$ in Marcos's view (Def. \ref{def:marcos}). 
\item $\e$  is paracomplete in  Omori's view (Def. \ref{def:omori}).
\end{enumerate}
\end{corollary}
\begin{proof} Let $W=\{w,s,t\}$ and $p$ be a propositional variable, i.e., $p\in \mathrm{Var}$. Then, suppose $\mm=\langle W,N,V \rangle$ where $N(w)=\{\{w,s\}, \{s,t\}\}$ and  
 $V$ is a valuation such that $V(p)=\{s,t\}$.  Then, we obtain  
    $\mm,w\not\models p$ , $\mm,w\not\models \nb p$ and $\mm,w\not\models \nd p$. \begin{enumerate}[label=(\roman*)]\item $\nb$ and $\nd$ qualify as  decent negations by Lemma \ref{lem:decent neg}. The model $\mm$ shows that $\e$ is $\nb$-incomplete and $\nd$-incomplete. So, $\not\models p \lor\nb p$ and  $\not\models p\lor \nd p$.  Also, $\mm, w\models \neg p$. This implies $\neg p\not\models p$. Thus,  $\e$ is non-overcomplete.  Hence, $\e$ is paracomplete with respect to negations $\nb$ and $\nd$ according to Definition   \ref{def:marcos} \ref{def:pcomplete Marcos}.

    \item   $\e$ extends $\mathbf{CPL}$ with the language $\lang'$.  The model  $\mm$ clearly shows that $\not\models p \lor\nb p$ and $\not\models p \lor \nd p $ where $p$ is a propositional variable. By soundness of $\e$ with respect to neighbourhood  semantics, $\e$ is paracomplete according to the Definition \ref{def:omori} \ref{def:pcomplete Omori}. \qed
    \end{enumerate}
    
\end{proof}

\begin{remark}
     $\e$ is a paranormal logic in view of Marcos by Definition \ref{def:marcos} \ref{def:pnormal Marcos} with respect to the negations $\nb$ and $\nd$. 
\end{remark}
Let $\mathscr{L}=\langle\lang, \models_\mathscr{L}\rangle$ be an arbitrary logic. Suppose $\mathscr{L}$ is equipped with a unary operator $\sim$. Then,  the following rule is called the \textit{contraposition rule} (or global contraposition rule) with respect to $\sim$.
\label{contra rule}For $\alpha,\beta\in \lang$,  $$\dfrac{\alpha\models_\mathscr{L} \beta}{\sim\beta\models_\mathscr{L} \sim\alpha}$$

\begin{theorem}\label{thm:contra rule failure in E}
The contraposition rule fails with respect to both $\nb$ and $\nd$ in $\e$, i.e., for $\sim\in\{\nb,\nd\}$, there exist formulas $\alpha,\beta\in \lang'$ such that $\alpha\models\beta$ but $\sim\beta\not\models\sim\alpha$.
    \end{theorem}
    \begin{proof}
        Let $p,q\in \mathrm{Var}$. 
\begin{enumerate}[label=(\roman*)]
\item Suppose  $W_1=\{w,s\}$ and  $\mm_1=\langle W_1,N_1,V_1\rangle$ where  $V_1$ is a valuation such that $V_1(p)=\{w,s\}$, $V_1(q)=\{w\}$ and $N_1(w)=\{\{w\}\}$. We know that $p\land q\models p$ in $\e$.  $\mm_1, w\models \nb p $ as  $\llbracket p\rrbracket_{\mm_1}\not \in N_1(w)$, and  $\mm_1,w\not\models \nb(p\land q)$ as $\llbracket p\land q\rrbracket_{\mm_1}\in N_1(w)$. Thus, $\nb p \not\models \nb(p \land q)$.
 \item Suppose  $W_2=\{w,s,t\}$ and let $\mm_2$ be a model, $\mm_2=\langle W_2,N_2,V_2\rangle$ where  $V_2$ is a valuation such that $V_2(p)=\{w\}$, $V_2(q)=\{s\}$ and $N_2(w)=\{\{s,t\}\}$. We know that $p\land q\models p$ in $\e$.  $\mm_2, w\models \nd p $ as  $\llbracket p\rrbracket^C_{\mm_2} \in N_2(w)$, and  $\mm_2,w\not\models \nd(p\land q)$ as $\llbracket p\land q \rrbracket^C_{\mm_2}\not \in N_2(w)$. Thus, $\nd p \not\models \nd(p \land q)$.\hfill \qed
\end{enumerate}
\end{proof}

Let $\mathscr{L}=\langle\lang, \models_\mathscr{L}\rangle$ be an arbitrary logic. Suppose $\mathscr{L}$ is equipped with a unary operator $\sim$. Then,  the following formulas are called De Morgan's laws with respect to $\sim$.
\label{De Morgan Laws} For $\alpha,\beta\in \lang$, 
\begin{enumerate}[label=(\roman*)]
    \item $\sim(\alpha\lor\beta)\leftrightarrow (\sim \alpha \land \sim \beta)$
    \item $\sim(\alpha\land\beta)\leftrightarrow (\sim \alpha \lor \sim \beta)$
\end{enumerate}

In the following theorem, we prove that these laws fail in $\e$. More specifically, we show that none of the implications holds.  

\begin{theorem}\label{thm:demorgan failure in E}
    For $\sim\in\{\nb, \nd\}$, there exist formulas $\alpha,\beta\in \lang'$ such that the following formulas are not valid in $\e$. 

   \noindent
        (i) \label{i demorgan} $\sim(\alpha\lor\beta)\rightarrow (\sim \alpha \land \sim \beta)$;\qquad
        (ii) \label{ii demorgan}$(\sim \alpha \land \sim \beta)\rightarrow \sim(\alpha\lor\beta) $;

    \noindent
         (iii) \label{iii demorgan} $\sim(\alpha\land\beta)\rightarrow (\sim \alpha \lor \sim \beta)$;\qquad
         (iv) \label{iv demorgan}$(\sim \alpha \lor \sim \beta)\rightarrow\sim(\alpha\land\beta) $.
    
\end{theorem}

\begin{proof}  Suppose that $\alpha\models\beta$, then for any neighbourhood  model $\mm$, $\llbracket\alpha\lor\beta\rrbracket_\mm=\llbracket\beta\rrbracket_\mm$ and $\llbracket\alpha\land\beta\rrbracket_\mm=\llbracket\alpha\rrbracket_\mm$. By semantic conditions for $\nb$ and $\nd$, we obtain,  $\sim(\alpha\lor\beta)\equiv \sim\beta$ and $\sim(\alpha\land\beta)\equiv \sim \alpha$ for $\sim\in\{\nb,\nd\}$.  If (i) is valid in $\e$, then, together with $\sim(\alpha\lor\beta)\equiv \sim\beta$, we obtain $\sim\beta\models\sim\alpha$ which is a contradiction to Theorem \ref{thm:contra rule failure in E}. Similarly, if we assume that (iv) is valid in $\e$, then with $\sim(\alpha\land\beta)\equiv \sim \alpha$, we obtain $\sim\beta\models\sim\alpha$, which is a contradiction to Theorem \ref{thm:contra rule failure in E}.  Hence, (i) and (iv) are not valid in $\e$.

Suppose $p,q\in \mathrm{Var}$ and let $W=\{w,s,t\}$. We omit the detailed explanation and simply provide the models. One can easily verify. 
 Let $N$ be a neighbourhood  function such that  $N(w)=\{\{w,t\}, \{w,s\}\}$, $N(s)=\{\{s\}\}$, and $N(t)=\{\{t\}\}$. Let $V_1$, $V_2$, $V_3$ and $V_4$ be valuations such that $V_1(p)=\{w\}$, $V_1(q)=\{s\}$, $V_2(p)=\{s\}$, $V_2(q)=\{t\}$, $V_3(p)=\{w,s\}$, $V_3(q)=\{w,t\}$, $V_4(p)=\{w,s\}$ and $V_4(q)=\{s,t\}$.
    
Let $\mm_1=\langle W,N,V_1 \rangle$. 
Then, $\mm_1,w \not\models \nb p\land \nb q \rightarrow \nb (p \lor q)$. 
Let  $\mm_2=\langle W,N,V_2 \rangle$. Then, 
     $\mm_2,w \not\models \nd p\land \nd q \rightarrow \nd (p \lor q)$. 
 Let $\mm_3=\langle W,N,V_3 \rangle$. Then, $\mm_3,w \not\models \nb( p\land q) \rightarrow (\nb p \lor \nb q)$. 
 Suppose $\mm_4=\langle W,N,V_4 \rangle$. Then  $\mm_4,w \not\models \nd( p\land q) \rightarrow (\nd p \lor \nd q)$. Hence, (ii) and (iii) are not valid in $\e$. 
    \qed
   \end{proof}

We observe that the negations $\nb$ and $\nd$ in $\e$ have weaker structural properties in the sense that they satisfy neither the contraposition rule nor any of De Morgan's laws. This is not the case for negations induced by modal operators in normal modal logics. This will be discussed later in Section \ref{sec:comparison}. 

\begin{remark}\label{rem:jimde groot marcos}
   As mentioned in \cite{MR2193566}, there is no unique answer in the literature to the question of what rules every negation should obey. So, Marcos introduced the concept of decent negation. In one of the recent works \cite{de2025intrinsic}, the authors focus on logics containing a unary operator $\sim$ and call it a negation if  contraposition holds, i.e., whenever one formula follows from another,  the $\sim$
attached to the latter follows from the $\sim$ attached to the former. In \cite{mruczek2008paraconsitent}, the authors use B\'{e}ziau's negation (possibly not $\Diamond\neg$) to construct paraconsistent logics expressible by  Kripke semantics where non-normal worlds are allowed.  However, they begin by defining a class of logics extending the positive fragment of $\mathbf{CPL}$  equipped with a negation satisfying the contraposition rule. In contrast to these, $\nb$ and $\nd$ do not satisfy the contraposition rule (see Theorem \ref{thm:contra rule failure in E}) in $\e$. So, in this article, we analyze negations $\nb$ and $\nd$ that do not obey the rule.  
\end{remark}
\section{LFIs and LFUs}\label{sec:lfi and lfu}
Logics of Formal Inconsistency (LFIs) and logics of formal undeterminedness (LFUs) were introduced in \cite{MR2008227}. LFIs are paraconsistent logics, i.e., logics with non-explosive negations. Moreover, these logics use a unary \emph{consistency} operator, usually denoted by $\circ$.
In this section, we prove that $\e$ gives rise to LFIs and LFUs with respect to both  the modal negations $\nb$ and $\nd$. Before proving the main results, we present a few definitions from \cite{MR2008227}.

\begin{definition}\label{dfn:Tarski}
A logic $\logl=\langle\lang,\bvdash_\logl\rangle$ is called \emph{Tarskian} if it satisfies the following properties. For every $\Gamma\cup\Delta\cup\{\alpha\}\subseteq\lang$, if $\alpha\in\Gamma$, then $\Gamma\bvdash_\logl\alpha$ (Reflexivity);
         if $\Gamma\bvdash_\logl\alpha$ and $\Gamma\subseteq\Delta$, then $\Delta\bvdash_\logl\alpha$ (Monotonicity);
         if $\Delta\bvdash_\logl\alpha$ and $\Gamma\bvdash_\logl\beta$ for all $\beta\in\Delta$, then $\Gamma\bvdash_\logl\alpha$ (Transitivity/Cut).
    
    $\logl$ is said to be \emph{finitary} if for every $\Gamma\cup\{\alpha\}\subseteq\lang$, $\Gamma\bvdash_\logl\alpha$ implies that there exists a finite $\Gamma_0\subseteq\Gamma$ such that $\Gamma_0\bvdash_\logl\alpha$.

    The logic $\logl$ is called \emph{standard} if it is Tarskian, finitary and structural.
\end{definition}
 \begin{remark}\label{rem:e std}
    By Definition \ref{def:syntactic consequence e}, it follows that $\e$ is Tarskian and finitary. By Remark \ref{rem:e structural}, $\e$ is structural, and hence standard.  
\end{remark}

\begin{definition}\label{def:lfi}
Let $\mathscr{L}=\langle\lang,\bvdash_\logl\rangle$ be a standard logic with a signature (a set of connectives) containing a negation $\sim$ and a primitive or defined unary consistency operator $\circ$. Then, $\mathscr{L}$ is said to be a \emph{Logic of Formal Inconsistency (LFI)} with respect to $\sim$ and $\circ$ if the following conditions hold.

\begin{enumerate}[label=(\roman*)]
    \item\label{lfi def i} $\{\alpha,\sim\alpha\}\not\bvdash_\logl\beta$ for some $\alpha,\beta\in\lang$.
    \item \label{lfi def ii}There exist $\alpha,\beta\in\lang$ such that \begin{enumerate}[label=(\alph*)]
        \item $\{\circ\alpha,\alpha\}\not\bvdash_\logl\beta$;
        \item $\{\circ\alpha,\sim\alpha\}\not\bvdash_\logl\beta$.
    \end{enumerate}
    \item \label{lfi def iii}$\{\circ\alpha,\alpha,\sim\alpha\}\bvdash_\logl\beta$ for all $\alpha,\beta\in\lang$.
\end{enumerate}
\end{definition}
\begin{definition}
    \label{def:lfu}
Let $\mathscr{L}=\langle\lang,\bvdash_\logl\rangle$ be a standard logic with a signature containing a negation $\sim$. Assume that $\mathscr{L}$ has a primitive or defined disjunction $\lor$ which enjoys the standard property, i.e., for every set of formulas $\Gamma\cup \{\alpha,\beta\}$, $Cn(\Gamma\cup \{\alpha\})\cap Cn(\Gamma\cup \{\beta\})=Cn(\Gamma\cup \{\alpha\lor\beta\})$ where $Cn(\Delta)=\{\gamma:\Delta\bvdash_\logl\gamma\}$, for every $\Delta$. Let $\hstar$ be a primitive or defined unary connective in the signature. Then, $\mathscr{L}$ is said to be a \emph{Logic of Formal Undeterminedness (LFU)} with respect to $\sim$ and $\hstar$ if the following conditions hold.

\begin{enumerate}[label=(\roman*)]
    \item \label{lfu dfn i}$\not\bvdash_\logl\alpha\lor\sim\alpha$ for some $\alpha\in\lang$.
    \item \label{lfu dfn ii}There exists $\alpha\in\lang$ such that \begin{enumerate}[label=(\alph*)]
        \item $\hstar\alpha\not\bvdash_\logl\alpha$;
        \item $\hstar\alpha\not\bvdash_\logl\sim\alpha$.
    \end{enumerate}
    \item\label{lfu dfn iii} $\hstar\alpha\bvdash_\logl\alpha\lor\sim\alpha$ for every $\alpha\in\lang$.
\end{enumerate}
\end{definition}

To show that $\e$ is an LFI with respect to $\nb$ and $\nd$, we need consistency operators that help recover ECQ with these negations in a non-trivial way ((ii) and (iii) of  Def. \ref{def:lfi}).
So, we define for any formula $\alpha\in \lang'$, $\circ_B \alpha:= \alpha \rightarrow \Box \alpha$ and $\circ_D \alpha:=\alpha\rightarrow\Diamond\alpha$.  
\begin{lemma}\label{lem:lfi nb nd condition}
There exist formulas $\alpha$, $\beta\in \mathcal{L}'$  such that the following hold in $\e$.
\begin{enumerate}[label=(\roman*)]
    \item \label{i nb lfi}$\cb\alpha, \alpha\not\bvdash\beta$;
    \item \label{ii nb lfi}$\cb\alpha,\nb\alpha\not\bvdash \beta$;
    \item \label{i nd lfi}$\cd\alpha, \alpha\not\bvdash\beta$;
    \item \label{ii nd lfi}$\cd\alpha,\nd\alpha\not\bvdash \beta$.
\end{enumerate}
\end{lemma}
\begin{proof}
    Let $p,q\in \mathrm{Var}$. Suppose $W=\{w,s\}$, $N(s)=\{\{s\}\}$ and $N(w)=\{\{w\}\}$. Let $V_1$ and $V_2$ be two valuations such that  $V_1(p)=\{w\}$,  $V_1(q)=\{s\}$, $V_2(p)=\{s\}$ and $V_2(q)=\{s\}$.
    
         Let $\mm=\langle W, N, V_1\rangle$.   Then, $\mm,w\models p$,  $\mm, w\models\Box p$ as $\llbracket p\rrbracket_{\mm}\in N(w)$. This implies $\mm, w\models \cb p$.
        Also, $\mm,w\models\Diamond p$ as $\llbracket p\rrbracket^C_{\mm}\not \in N(w)$. This implies $\mm, w\models \cd p$.
        But $\mm, w\not\models q$ as $ w\not \in V(q)$.   Consequently, we obtain $\cb p, p \not\models q$ and $\cd p, p \not\models q$.  Soundness of $\e$ gives $\cb p, p \not\bvdash q$ and $\cd p, p \not\bvdash q$ thereby proving \ref{i nb lfi} and \ref{i nd lfi}.
        
             Let $\mm'$ be a model, $\mm'=\langle W, N, V_2\rangle$.  Then, $\mm', w\not\models p$  implies $\mm', w\models \cb p$ and $\mm', w\models \cd p$. $\mm', w\models \nb p$ as $V_2(p)\not\in N(w)$ and $\mm', w\models\nd p$ as  $(V_2(p))^C\in N(w)$. Also,  $\mm', w\not\models q$ as $w\not \in V_2(q)$.  Hence, we get  $\cb p, \nb p \not\models q$ and  $\cd p, \nd p \not\models q$.  Then, by soundness of $\e$, \ref{ii nb lfi} and \ref{ii nd lfi} are proved.
             \qed
           \end{proof}

\begin{theorem}
\begin{enumerate}[label=(\roman*)]
\item  $\e$ is an LFI with respect to $\nb$ and $\cb$.
   \item  $\e$ is an LFI with respect to $\nd$ and $\cd$.
    \end{enumerate}
\end{theorem}
\begin{proof} By Remark \ref{rem:e std}, $\e$ is standard. ECQ fails with respect to $\nb$ in $\e$ by Theorem \ref{thm:ecq fails}.  This proves \ref{lfi def i} of Def. \ref{def:lfi}. One can directly derive $\{\cb\alpha,\alpha,\nb\alpha\}\bvdash \beta$ for any $\alpha,\beta\in \lang'$ using ECQ (for $\neg$)   and the rule (MP) in $\e$. Thus, we get \ref{lfi def iii} of Def. \ref{def:lfi}. Lemma \ref{lem:lfi nb nd condition} (\ref{i nb lfi} and \ref{ii nb lfi}) give \ref{lfi def ii} of Def. \ref{def:lfi}. Thus, $\e$ is an LFI with respect to $\nb$ and $\cb$. 

     With similar arguments, one can prove that $\e$ is an LFI with respect to $\nd$ and $\cd$ using Theorem \ref{thm:ecq fails}  and Lemma \ref{lem:lfi nb nd condition}.  
     \qed
\end{proof}

This shows that $\e$ provides two different LFIs. Now, we show that $\e$ is an LFU with respect to both negations $\nb$ and $\nd$. We define determinedness operators, that recover LEM in a non-trivial way ((ii) and (iii) of Definition \ref{def:lfu}). For any formula $\alpha\in \lang'$, we define $\db\alpha:=\Box \alpha\rightarrow\alpha$ and $\dD\alpha:=\Diamond \alpha\rightarrow  \alpha$.   
\begin{lemma}\label{lem:lfu nb nd condition}
    There exists a formula $\alpha\in \mathcal{L}'$ such that the following hold in $\e$.
\begin{enumerate}[label=(\roman*)]
    \item \label{i nb lfu}$\db\alpha\not\bvdash\alpha$
    \item \label{ii nb lfu}$\db\alpha\not\bvdash \nb \alpha$ 
    \item\label{i nd lfu} $\dD\alpha\not\bvdash\alpha$
    \item\label{ii nd lfu} $\dD\alpha\not\bvdash \nd \alpha$
\end{enumerate}
\end{lemma}
\begin{proof} Let $p \in \mathrm{Var}$. Suppose $W=\{w,s\}$ and $N(w)=\{\{w\}\}$. Suppose $V_1$ and $V_2$ are two valuations such that $V_1(p)=\{s\}$ and $V_2(p)=\{w\}$. 

    Suppose $\mm=\langle W, N, V_1\rangle$. Then, $\mm, w\not\models\Box p$. This implies $\mm,w\models \db p$. Moreover,  $\mm, w\not\models\Diamond p$ as $(V_1(p))^C\in N(w)$. Hence $\mm, w\models \dD p$. But $\mm,w\not\models p$. Therefore, $\db p\not\models p$ and $\dD p\not\models p$.  By soundness of $\e$, \ref{i nb lfu} and \ref{i nd lfu} are proved. 

  Let $\mm'$ be a model, $\mm'=\langle W, N, V_2\rangle$. Then $\mm',w\models p$. This implies  $\mm',w\models \db p$ and $\mm',w\models \dD p$. But $\mm',w\not\models\nb p$.  Since $\llbracket\neg p\rrbracket_{\mm'}\not\in N(w)$, we get $\mm',w\not\models\nd p$. Therefore, $\db p\not\models\nb p$ and $\dD p\not\models \nd p$. By soundness of $\e$,  \ref{ii nb lfu} and \ref{ii nd lfu} are proved.
  \qed
  \end{proof}

\begin{theorem}
\begin{enumerate}[label=(\roman*)]
\item  $\e$ is an LFU with respect to $\nb$ and $\db$.
   \item  $\e$ is an LFU with respect to $\nd$ and $\dD$.
    \end{enumerate}
\end{theorem}
\begin{proof}
 By Remark \ref{rem:e std}, $\e$ is a standard logic. $\e$ extends $\mathbf{CPL}$ and hence $\lor$ behaves classically. The property can be verified directly using neighbourhood semantics as $\e$ is sound and complete with respect to neighbourhood semantics. So, $\e$ enjoys the standard property of $\lor$.
 
LEM fails in $\e$ with respect to $\nb$ by Theorem \ref{thm:lem fails}. This proves condition \ref{lfu dfn i} of Def. \ref{def:lfu}. Lemma \ref{lem:lfu nb nd condition} proves condition \ref{lfu dfn ii} of Def. \ref{def:lfu}.   For any formula $\alpha\in \mathcal{L'}$, $\db \alpha=\Box\alpha\rightarrow \alpha$.  It is known that $\Box\alpha\rightarrow \alpha\equiv\neg\Box \alpha\lor \alpha$. This implies $\db \alpha\models \neg\Box \alpha\lor \alpha $. By completeness of $\e$, $\db\alpha\bvdash \nb \alpha\lor \alpha$. This is true for all formulas $\alpha\in \mathcal{L'}$. This proves condition \ref{lfu dfn iii} of Def. \ref{def:lfu}. So, $\e$ is an LFU with respect to $\nb$ and $\db$.
 
     Similarly,  $\e$ is an LFU with respect to $\nd$ and $\dD$ using Theorem \ref{thm:lem fails} and Lemma \ref{lem:lfu nb nd condition}.\hfill \qed
    
    \end{proof}
\begin{remark}\label{rem:Rmbc}
In \cite{MR4469283}, the authors presented an LFI $\mathbf{RmbC}$ and proved that $\mathbf{RmbC}$ can be defined within the bimodal non-normal logic $\e\oplus\e$, the fusion of the non-normal modal logic $\e$ with itself.  $\e\oplus\e$ has two independent modalities $\Box_1$ and $\Box_2$ required for defining the paraconsistent negation and consistency connective. An extension of  $\mathbf{RmbC}$ called $\mathbf{RmbCciw}$, in which consistency coincides with non-contradictoriness, is definable within $\e$. It would be interesting to find connections between  $\mathbf{RmbCciw}$ and the LFIs obtained here from $\e$. We leave this for future work.  
However, their modal representation interprets the paraconsistent negation $\neg_{\mathbf{RmbC}}$ of a formula $\alpha$ by $\alpha\rightarrow\Box_1\alpha$. In this article, we study the negations $\nb$ and $\nd$ induced directly by the modal operator of $\e$.
\end{remark}
\section{Comparison with other logics}\label{sec:comparison}
 
 Some formulas and rules listed in the beginning of Section \ref{sec:system e} are not valid in $\e$. Such countermodels can be found in  \cite[Chapter 2]{pacuit2017neighborhood}. All of them involve modal operators, so they also affect modal negations $\nb$ and $\nd$. In this section, we discuss how the behaviour of the modal negations changes if we add some of these formulas to $\e$ and compare $\e$ with several well-known logics like Ja\'{s}kowski's Discursive logic $D_2$, intuitionistic logic, and minimal logic.  

\vspace{0.75mm}

We begin by adding the axiom scheme (M) \ref{M axiom} to  $\e$. The resultant system is called  $\mathbf{EM}$. It extends the system $\e$ by adding all instances of the axiom scheme (M). $\mathbf{EM}$ is sound and complete with respect to the class of monotonic\footnote{ Let $W$ be any non-empty set. Then, a collection $\mathcal{U}$ of subsets  of  $W$ is called monotonic if this condition holds: for each $X\subseteq W$, if $X\in \mathcal{U}$, and $X\subseteq Y\subseteq W$, then  $Y\in \mathcal
{U}$. A frame $\langle W, N\rangle$ is called monotonic if for each $w\in W$, $N(w)$ is monotonic. A proof of completeness theorem can be found in Chapter 2 of \cite{pacuit2017neighborhood}.} frames. Thus, if $\alpha\models_{\mathbf{EM}}\beta$, we obtain $\Box\alpha\models_{\mathbf{EM}}\Box\beta$, as is easily verified. By the contraposition principle of classical negation, 
$\neg \Box\beta\models_{\mathbf{EM}}\neg \Box\alpha$ or $\nb\beta\models_{\mathbf{EM}} \nb\alpha$. Similarly, we get  $\nd\beta\models_{\mathbf{EM}}\nd\alpha$ if $\alpha\models_{\mathbf{EM}}\beta$.

Thus, $\mathbf{EM}$ validates the contraposition rule for both $\nb$ and $\nd$,  whereas this rule fails in $\e$ (Theorem \ref{thm:contra rule failure in E}). Moreover, $\mathbf{EM}$ validates  one direction of each De Morgan's law for $\nb$ and $\nd$, namely $\sim\alpha\lor \sim \beta \rightarrow \sim (\alpha\land \beta)$ and $\sim (\alpha\lor \beta)\rightarrow (\sim\alpha\land \sim \beta)$ where $\sim\in\{\nb,\nd\}$ due to (M).

Adding the axiom scheme (C) \ref{C axiom} to $\mathbf{EM}$ yields $\mathbf{EMC}$.  Using $\neg$ contraposition and De Morgan's laws on (C), we obtain $\models_{\mathbf{EMC}}\nb(\alpha\land\beta)\rightarrow (\nb\alpha\lor \nb \beta)$. The converse is already valid because of (M). Thus, $\models_{\mathbf{EMC}}\nb(\alpha\land\beta)\leftrightarrow (\nb\alpha\lor \nb \beta)$. Similarly, we obtain $\models_{\mathbf{EMC}}\nd(\alpha\lor \beta)\leftrightarrow \nd\alpha\land\nd\beta$.  Thus, one De Morgan's law is valid in the system $\mathbf{EMC}$ with respect to the modal negations $\nb$ and $\nd$.

  (K) is valid in $\mathbf{EMC}$ but $\mathbf{EMC}$  does not coincide with $\mathbf{K}$ because the rule (Nec) is not valid in $\mathbf{EMC}$. Note that  $\mathbf{K}$ does not restore ECQ and LEM with respect to $\nb$ and $\nd$. But since $\nb$ and $\nd$ in $\mathbf{EMC}$ and $\mathbf{K}$ satisfy contraposition, systems of this kind have already been studied in the literature, as mentioned in Remark \ref{rem:jimde groot marcos}.

Ja\'{s}kowski's Discussive logic $\mathbf{D_2}$ \cite{jaskowski1969propositional} is a very well-known paraconsistent logic. The logic $\mathbf{S5}$ is used to establish paraconsistency in $\mathbf{D_2}$. The main idea is as follows. A conclusion $q$ discussively follows from formulas $p_1,  p_2, \dots, p_n$ if $\Diamond q$ follows from $\Diamond p_1$, $\Diamond p_2$,\dots, $\Diamond p_n$ in $\mathbf{S5}$. As a result, there exist formulas $\alpha$, $\beta$ in the language of $\mathbf{D_2}$, such that $\alpha,\neg\alpha\not\bvdash_\mathbf{D_2} \beta$. We observe that paraconsistency in  $\mathbf{D_2}$ is due not only to weak negation but also to the defined consequence relation. Therefore, $\alpha\land\neg\alpha\bvdash_\mathbf{D_2}\beta$ for all formulas $\alpha$ and $\beta$, despite the fact that $\alpha,\neg \alpha\bvdash_\mathbf{D_2}\beta$ does not hold.  However, in $\e$, paraconsistency arises due to the behaviour of negations $\nb$ and $\nd$.  

Intuitionistic logic was developed by Brouwer and formally axiomatized by Heyting. It is one of the famous paracomplete logics. It provides a logical basis for constructive mathematics. The formal system of intuitionistic propositional logic is usually denoted by $\mathbf{IPC}$. The signature used for $\mathbf{IPC}$ is $\{\lor, \land, \rightarrow, \bot \}$, where $\bot$  is a primitive constant. Negation $\neg$ of a formula $\alpha$ is defined as $\neg\alpha:=\alpha\rightarrow \bot$. It is well known that  in $\mathbf{IPC}$, LEM fails, i.e., $\not\bvdash_{\mathbf{IPC}} \alpha\lor \neg\alpha$  for some formula $\alpha$ and  ECQ holds, i.e., $\alpha, \neg\alpha \bvdash_{\mathbf{IPC}}\beta$ for all formulas $\alpha$ and $\beta$ in the language. Thus, it is paracomplete but not paraconsistent with respect to $\neg$.  

\noindent
 Minimal logic ($\mathbf{ML}$) was introduced by Johansson in \cite{Johansson1937}.  $\mathbf{ML}$ has the same language as $\mathbf{IPC}$. Both ECQ and LEM fail in $\mathbf{ML}$, i.e.,  $\alpha,\neg \alpha\not\bvdash_{\mathbf{ML}} \beta$ for some formulas $\alpha$ and $\beta$, and   $\not\bvdash_{\mathbf{ML}}\alpha\lor\neg\alpha$ for some $\alpha$.
 
 But the negations of $\mathbf{IPC}$ and $\mathbf{ML}$ have more properties than $\nb$ and $\nd$ of $\e$. The principle of contraposition and one De Morgan's law ($\neg(\alpha\lor\beta)\leftrightarrow \neg \alpha\land \neg \beta$) hold in both $\mathbf{IPC}$ and $\mathbf{ML}$. In contrast, $\e$ does not satisfy any of De Morgan's laws (Theorem \ref{thm:demorgan failure in E}) and the contraposition rule (Theorem \ref{thm:contra rule failure in E}). For further details about $\mathbf{IPC}$ and $\mathbf{ML}$, one can see \cite{bezhanishvili2006intuitionistic,troelstra2000basic}.

\paragraph{\textbf{Conclusion.}} 
We observe that $\e$ yields two LFIs and two LFUs with respect to the modal negations $\nb$ and $\nd$, together with suitable consistency operators $\cb$ and $\cd$, and corresponding determinedness operators $\db$ and $\dD$.  We note that even weak non-normal modal systems can qualify as LFIs and LFUs. One can explore the connections between $\e$ and $\mathbf{RmbCciw}$, as discussed in Remark \ref{rem:Rmbc}. A natural direction for future work is to investigate how these weak, non-contrapositive modal negations behave in dynamic or multi-agent non-normal modal frameworks where agents are not logically omniscient. 
\begin{credits}
\subsubsection{\ackname} The author thanks her advisor, Sankha S. Basu, for his helpful comments and suggestions on earlier drafts of this article. In particular, general discussions on the theory of negation provided important background for this work. The author also thanks the anonymous reviewers for their helpful feedback.

\end{credits}

\bibliographystyle{splncs04}
\bibliography{fntp}

@article {MR2193566,
    AUTHOR = {Marcos, {J}o\~{A}o},
     TITLE = {Nearly every normal modal logic is paranormal},
   JOURNAL = {Logique et Anal. (N.S.)},
  FJOURNAL = {Logique et Analyse. Nouvelle S\'erie},
    VOLUME = {48},
      YEAR = {2005},
    NUMBER = {189-192},
     PAGES = {279--300},
      ISSN = {0024-5836},
   MRCLASS = {03B45},
  MRNUMBER = {2193566},
MRREVIEWER = {Branislav\ Bori\v ci\'c},
}

@book{pacuit2017neighborhood,
  title= {Neighbourhood semantics for modal logic},
  author={Pacuit, Eric},
  year={2017},
  publisher={Springer},
  doi={10.1007/978-3-319-67149-9}
}

@article{homori2025s5,
  title={Is {S}5 Paraconsistent?},
  author={Omori, Hitoshi},
  journal={Principia: an international journal of epistemology},
  volume={29},
  number={2},
  pages={243-254},
  year={2025},
  publisher={Universidade Federal de Santa Catarina (UFSC)},
  doi={10.5007/1808-1711.2025.e107472}
}

@article{mruczek2008paraconsitent,
  title={Paraconsitent logics obtained by {J.}-{Y.} {B}{\'e}ziau’s method by means of some non-normal modal logics},
  author={Mruczek-Nasieniewska, Krystyna and Nasieniewski, Marek},
  journal={Bulletin of the Section of Logic},
  volume={37},
  number={3/4},
  pages={185--196},
  year={2008}
}

@article{beziau2002s5,
  title={S5 is a paraconsistent logic and so is first-order classical logic},
  author={B{\'e}ziau, J-Y},
  journal={Logical Studies},
  volume={9},
  pages={301--309},
  year={2002},
  publisher={Федеральное государственное бюджетное учреждение науки {\guillemotleft}Институт философии~…}
}

@article{de2025intrinsic,
  title={Intrinsic and relative characterization results for logics with negative modalities},
  author={de Groot, Jim and Marcos, Jo{\~a}o and Stefanes, Rodrigo},
  journal={arXiv preprint arXiv:2512.15496},
  year={2025},
  doi={10.48550/arXiv.2512.15496}
}

@article {MR4939722,
    AUTHOR = {Basu, Sankha S. and Roy, Sayantan},
     TITLE = {Generalized explosion principles},
   JOURNAL = {Studia Logica},
  FJOURNAL = {Studia Logica. An International Journal for Symbolic Logic},
    VOLUME = {113},
      YEAR = {2025},
    NUMBER = {4},
     PAGES = {1025--1060},
      ISSN = {0039-3215,1572-8730},
   MRCLASS = {03B53},
  MRNUMBER = {4939722},
       DOI = {10.1007/s11225-024-10153-x},
}

@article{jaskowski1969propositional,
  title={Propositional calculus for contradictory deductive systems (communicated at the meeting of march 19, 1948)},
  author={Ja{\'s}kowski, Stanis{\l}aw},
  journal={Studia Logica: An International Journal for Symbolic Logic},
  volume={24},
  pages={143--160},
  year={1969},
  publisher={JSTOR},
  doi={10.1007/BF02134311}
}

@book{troelstra2000basic,
  title={Basic proof theory},
  author={Troelstra, Anne Sjerp and Schwichtenberg, Helmut},
  edition   = {2},
  series    = {Cambridge Tracts in Theoretical Computer Science},
  number    = {43},
  publisher = {Cambridge University Press},
  address   = {Cambridge},
  year      = {2000},
  doi       = {10.1017/CBO9781139168717}
  }

@techreport{bezhanishvili2006intuitionistic,
  title={Intuitionistic logic},
  author={Bezhanishvili, Nick and De Jongh, Dick},
   institution = {Institute for Logic, Language and Computation, Universiteit van Amsterdam},
  type        = {Prepublication Series},
  number      = {PP-2006-25},
  year={2006}
}

@article{jaskowski1948rachunek,
  title={Rachunek zdan dla system{\'o}w dedukcyjnych sprzecznych},
  author={Ja{\'s}kowski, Stanis{\l}aw},
  journal={Studia Societatis Scientiarum Torunensis, Sectio A},
  volume={1},
  number={5},
  pages={57--77},
  year={1948}
}

@incollection{scott1970advice,
  title={Advice on modal logic},
  author={Scott, Dana},
  booktitle={Philosophical problems in logic: Some recent developments},
  pages={143--173},
  year={1970},
  publisher={Springer},
  doi={10.1007/978-94-010-3272-8_7}
}

@article{montague1970universal,
  author  = {Richard Montague},
  title   = {Universal Grammar},
  journal = {Theoria},
  volume  = {36},
  number  = {3},
  pages   = {373--398},
  year    = {1970},
  doi     = {10.1111/j.1755-2567.1970.tb00434.x}
}

@article {MR4469283,
    AUTHOR = {Carnielli, Walter and Coniglio, Marcelo E. and Fuenmayor,
              David},
     TITLE = {Logics of formal inconsistency enriched with replacement: an
              algebraic and modal account},
   JOURNAL = {Rev. Symb. Log.},
  FJOURNAL = {The Review of Symbolic Logic},
    VOLUME = {15},
      YEAR = {2022},
    NUMBER = {3},
     PAGES = {771--806},
      ISSN = {1755-0203,1755-0211},
   MRCLASS = {03B45 (03B62 03G27)},
  MRNUMBER = {4469283},
       DOI = {10.1017/S1755020321000277},
      
}

@incollection {MR2008227,
    AUTHOR = {Carnielli, Walter A. and Marcos, {J}o\~{A}o},
     TITLE = {A taxonomy of {${\bf C}$}-systems},
 BOOKTITLE = {Paraconsistency ({S}\~ao {S}ebasti\~ao, 2000)},
    SERIES = {Lecture Notes in Pure and Appl. Math.},
    VOLUME = {228},
     PAGES = {1--94},
 PUBLISHER = {Dekker, New York},
      YEAR = {2002},
      ISBN = {0-8247-0805-9},
   MRCLASS = {03B53 (03B45 03B50)},
  MRNUMBER = {2008227},
MRREVIEWER = {Jean-Yves\ B\'eziau},
       DOI = {10.1201/9780203910139-3},
      
}

@book{chellas1980modal,
  title={Modal logic: an introduction},
  author={Chellas, Brian F},
  volume={316},
  year={1980},
  publisher={Cambridge university press Cambridge},
  doi={10.1017/CBO9780511621192}
}

@article{holliday2015epistemic,
  title={Epistemic closure and epistemic logic {I}: Relevant alternatives and subjunctivism},
  author={Holliday, Wesley H},
  journal={Journal of Philosophical Logic},
  volume={44},
  number={1},
  pages={1--62},
  year={2015},
  publisher={Springer},
  doi={10.1007/s10992-013-9306-2}
}

@article{dosen1984negative,
  title={Negative modal operators in intuitionistic logic},
  author={Do\v{s}en, Kosta},
  journal={Publ. Inst. Math., Nouv. Ser},
  volume={35},
  number={49},
  pages={15--20},
  year={1984}
}

@article{hansen2009neighbourhood,
  title={Neighbourhood structures: Bisimilarity and basic model theory},
  author={Hansen, Helle Hvid and Kupke, Clemens and Pacuit, Eric},
  journal={Logical Methods in Computer Science},
  volume={5},
  year={2009},
  publisher={Episciences. org},
  doi={10.2168/LMCS-5(2:2)2009}
}

@book{segerberg1971essay,
  title={An essay in classical modal logic},
  author={Segerberg, Karl Krister},
  year={1971},
  publisher={Stanford University}
}

@article{Johansson1937,
  author  = {Johansson, Ingebrigt},
  title   = {Der {Minimalkalk\"ul}, ein reduzierter
             intuitionistischer {Formalismus}},
  journal = {Compositio Mathematica},
  volume  = {4},
  pages   = {119--136},
  year    = {1937}
}

@inproceedings{Pnueli1977,
  author    = {Pnueli, Amir},
  title     = {The Temporal Logic of Programs},
  booktitle = {18th Annual Symposium on Foundations of Computer Science (FOCS)},
  year      = {1977},
  pages     = {46--57},
  publisher = {IEEE},
  doi       = {10.1109/SFCS.1977.32}
}

@inproceedings{RaoGeorgeff1995,
  author    = {Rao, Anand S. and Georgeff, Michael P.},
  title     = {{BDI} Agents: From Theory to Practice},
  booktitle = {Proceedings of the First International Conference on Multi-Agent Systems (ICMAS)},
  year      = {1995},
  pages     = {312--319},
  publisher = {AAAI Press}
}

@article{vonWright1951,
  author  = {von Wright, Georg Henrik},
  title   = {Deontic Logic},
  journal = {Mind},
  volume  = {60},
  number  = {237},
  pages   = {1--15},
  year    = {1951},
  doi     = {10.1093/mind/LX.237.1}
}
\end{document}